\documentclass[12pt]{article}

\usepackage{amsmath}
\usepackage{amsfonts}
\usepackage{amssymb}
\usepackage{eucal}

\newcommand{\disp}{\displaystyle}

\makeatletter

\@addtoreset{equation}{section}
\makeatother

\newtheorem{theorem}{Theorem}[section]

\newtheorem{lemma}[theorem]{Lemma}
\newtheorem{corollary}[theorem]{Corollary}

\newtheorem{remark}[theorem]{Remark}

\newcommand{\qed}{\hbox{\rule{6pt}{6pt}}}

\begin{document}

\title{
Remarks on cutoff phenomena for random walks on Hamming Schemes}
\vspace{1pc}
\author{Katsuhiko Kikuchi}
\date{}
\maketitle

\footnote[0]{2010 {\it{Mathematics Subject Classification}}. 
Primary 05C81; Secondary 60C05, 05E18.}


\begin{abstract}
The sequence of the simple random walks on Hamming schemes 
$\{H(n, q)\}_{n=1}^{\infty}$ has a cutoff phenomenon 
for each integer $q$ greater than or equal to $3$. 
In this paper, for the sequence of simple random walks 
on Hamming schemes $\{H(n, q)\}_{n=1}^{\infty}$ with $q\geq 3$, 
we give a simple majorant and a sharp minorant function 
for total variance distances between transition distributions 
and stationary distributions. 

\end{abstract}

\section{Introduction.}
\quad\,\,
For many random walks on finite graphs, 
the transition distributions converge 
to the distributions of the equilibrium. 
Moreover, if we have useful majorant and minorant functions 
for the distance of them, we find the critical behavior 
of transition distributions, for example, the rapidity of decrease of the distances 
in the small range near the suitable time. 
Such the phenomenon and the time are called the cutoff phenomenon and 
the time to stationarity, respectively. 
In this paper, we give majorant and minorant functions for total variance distances 
between transition distributions and the distributions of the equilibrium 
for the simple random walks on Hamming schemes $\{H(n, q)\}_{n=1}^{\infty}$ 
with $q\geq 3$. 

Cutoff phenomenon is defined as follows. 
For a finite set $X$, we denote by $M(X)$ the vector space 
of all complex-valued measures on $X$. 
Take two measures $\mu , \nu \in M(X)$ on $X$ and 
define the total varialce distance $\Vert \mu -\nu \Vert _{TV}$ by 
\begin{equation*}
\Vert \mu -\nu \Vert _{TV}=\max \{|\mu (S)-\nu (S)|\,;\,S\subset X\}. 
\end{equation*}
Let $(X, E_X)$ be a simple connected finite unordered graph without loops. 
For $x, x'\in X$, we say that $x$ is adjacent to $x'$ if 
the (unordered) pair $\{x, x'\}$ belongs to $E_X$ and write $x\sim x'$. 
The transition probability $p(\cdot , \cdot )$ is a function on $X\times X$ 
such that (i)\,\,$p(x, x')\geq 0$, (ii)\,\,$p(x, x')>0$ if and only if $x\sim x'$, 
and (iii)\,\,$\disp{\sum _{x'\in X}p(x, x')=1}$ for any $x\in X$. 
For a nonnegative integer, we define transition probability $p^{(k)}(\cdot , \cdot )$ 
after $k$-steps recursively, by 
\begin{equation*}
p^{(0)}(x, x')=\delta _{x, x'},  \quad 
p^{(k)}(x, x')=\sum  _{y\in X}p^{(k-1)}(x, y)p(y, x'),\,\,k\geq 1, 
\end{equation*}
where $\delta _{x, x'}$ is the Kronecker delta. 
Fix an element $x^{(0)}\in X$ of $X$ and put 
$\nu ^{*k}(\cdot )=p^{(k)}(x^{(0)}, \cdot )$.  
Then, we see that $\nu ^{*k}$ is a probability measure on $X$. 
We say that the transition probability $p(\cdot , \cdot )$ is ergodic 
if there exists an integer $k_0$ such that $p^{(k)}(x, x')>0$ 
for any elements $x, x'\in X$ and $k\geq k_0$. 
A probability measure $\pi$ on $X$ is stationary if 
$\disp{\sum _{x'\in X}\pi (x')p(x', x)=\pi (x)}$ for any $x\in X$. 
Let $\{(X_n, E_{X_n})\}$ be a sequence of simple connected finite unordered graphs 
without loops, $\{p_n(\cdot , \cdot )\}$ the transition probabilities on $X_n$, 
$\{x_n^{(0)}\}$ the fixted points and $\{\pi _n\}$ the stationary probabilities. 
For sequences $\{a_n\}$, $\{b_n\}$ of positive real numbers 
with $\disp{\lim _{n\to \infty}\frac{b_n}{a_n}=0}$, 
the sequence of the Markov chains $\{(X_n, E_{X_n})\}$ has 
an $(a_n, b_n)$-cutoff if there exist functions $f_{\pm}\,:\,[0, +\infty )
\longrightarrow {\mathbb{R}}$ with $\disp{\lim _{c\to +\infty}f_+(c)=0}$, 
$\disp{\lim _{c\to +\infty}f_-(c)=1}$, and for each $c>0$, we have 
\begin{align}
\limsup _{n\to \infty}\Vert \nu _n^{*\lceil a_n+cb_n\rceil}-\pi _n\Vert _{TV} & 
\leq f_+(c), \label{0101} \\ 
\liminf _{n\to \infty}\Vert \nu _n^{*\lfloor a_n-cb_n\rfloor}-\pi _n\Vert _{TV} & 
\geq f_-(c), \label{0102}
\end{align}
where $\lceil \alpha \rceil$, $\lfloor \alpha \rfloor$ denote 
the least integer greater than or equal to $\alpha$, 
the greatest integer less than or equal to $\alpha$, respectively 
(see \cite{D1}, \cite{D2}, \cite{DS}). 
$f_+$ and $f_-$ are called the upper bound and the lower bound, and we often take 
as a monotone decreasing, a monotone increasing function, respectively. 
We remark that the existence of $f_+$ does not imply the ergodicity of each random walk 
on the graph $(X_n, E_{X_n})$ (see \cite{H}). 
If we take a majorant function $h_+$ for $\Vert \nu _n^{*\lfloor a_n+cb_n\rfloor}
-\pi _n\Vert _{TV}$ with $n\geq n_0$ for some positive integer $n_0$, 
we have the ergodicity of the random walk on $(X_n, E_{X_n})$ 
for any $n$ with $n\geq n_0$. 

Let $n$ be a positive integer and $q$ an integer with $q\geq 2$. 
we denote by $[q]_0=\{0, 1, \dots , q-1\}$ and $[n]=\{1, 2, \dots , n\}$. 
As a graph, the Hamming scheme $H(n, q)=(X_n, E_{X_n})$ is a finite graph 
with the vertex set 
$X_n=[q]_0^n$ and the edge set 
\begin{equation*}
E_{X_n}=\{\{x, x'\}\subset X_n\,;\,
\sharp \{j\in [n]\,;\,x_j\not= x'_j\}=1\}, 
\end{equation*}
where $x=(x_1, \dots , x_n)$, $x'=(x'_1, \dots , x'_n)\in X_n$, and 
$\sharp S$ is the cardinal number for a finite set $S$. 
The transition probability $p_n(\cdot , \cdot )$ 
for $(X_n, E_{X_n})$ is defined by 
\begin{equation*}
p_n(x, x')=\begin{cases}
\disp{\frac{1}{n(q-1)}}, & x\sim x', \\
\quad\,\,\,\, 0,         & {\text{otherwise}}. 
\end{cases}
\end{equation*}
Put $x^{(0)}=(0, \dots , 0)\in H(n, q)$ and 
$\nu _n(\cdot )=p_n(x^{(0)}, \cdot )$. 
Then, $\nu _n$ is a probability measure on $X_n$. 
The graph $H(n, q)$ is ergodic if and only if $q\geq 3$. 
Hora gives in \cite{H} the limit function $\disp{f_{\pm}(c)=Er\!f 
\left( \frac{e^{\mp \frac{c}{2}}}{2\sqrt{2}}\right)}$ 
of $\Vert \nu _n^{*(a_n\pm cb_n)}-\pi _n\Vert _{TV}$, 
where $Er\!f\,:\,{\mathbb{R}}\longrightarrow {\mathbb{R}}$ is 
the error function defined by $\disp{Er\!f(x)=\frac{2}{\sqrt{\pi}}
\int _0^xe^{-t^2}dt}$, and $\disp{(a_n, b_n)
=\!\left( \frac{n(q-1)}{2q}\log n(q-1), \frac{n(q-1)}{2q}\right)}$. 
Candidates of the majorant functions for $\Vert \nu _n^{*k}
-\pi _n\Vert _{TV}$ are given by Diaconis and Hanlon in \cite{DH}, 
and by Diaconis and Ram in \cite{DR}. 
In those papers, the simple random walk on a Hamming scheme 
$H(n, q)$ are regarded as a special case of Metropolis chains 
on a hypercube $H(n, 2)$. 
Mizukawa gives in \cite{M1} a majorant function for total variance distance 
with the different time to stationarity $\disp{\frac{n(q-1)}{2q}\log q^n}$ 
for the case $q\geq 3$, and in \cite{M2} a majorant
and a minorant function for the sequences of random walks 
on Hamming schemes with staying. 

For each integer $q$ with $q\geq 5$, we give a simple majorant function 
of $\Vert \nu _n^{*k}-\pi _n\Vert _{TV}$. 

\begin{theorem} \label{T0101}
Assume that $q\geq 5$. 
Let $\disp{k=\frac{n(q-1)}{2q}(\log n(q-1)+c)}$ be an integer with $c>0$. 
Then, we have 
\begin{equation} \label{0103}
\Vert \nu _n^{*k}-\pi _n\Vert _{TV}^2
\leq \frac{1}{4}(e^{e^{-c}}-1). 
\end{equation}
\end{theorem}

The key of the proof of this theorem is that $e^{-x}\geq |1-x|$ 
for any real number $x$ with $\disp{x\leq \frac{5}{4}}$. 
We cannot adapt the proof  for the case $q=3, 4$. 
The obstruction is that $|1-x|\geq e^{-x}$ 
for a real number $x$ with $\disp{x\geq \frac{4}{3}}$. 
So we replace the majorant function. 

\begin{theorem} \label{T0102}
{\rm{(1)}}\,\,
Assume that $q=3$ and $n\geq 3$. 
For any integer $\disp{k=\frac{n}{3}(\log 2n+c)}$ with $c>0$, we have 
\begin{equation} \label{0104}
\Vert \nu _n^{*k}-\pi _n\Vert _{TV}^2
\leq \frac{5}{2}(e^{e^{-c}}-1). 
\end{equation}
{\rm{(2)}}\,\,
Suppose that $q=4$ and $n\geq 2$. 
For each integer $\disp{k=\frac{3n}{8}(\log 3n+c)}$ with $c>0$, one has 
\begin{equation} \label{0105}
\Vert \nu _n^{*k}-\pi _n\Vert _{TV}^2
\leq \frac{9}{4}(e^{e^{-c}}-1). 
\end{equation}
\end{theorem}

While, we give minorant function for $\Vert \nu _n^{*k}
-\pi _n\Vert _{TV}$ as follows. 

\begin{theorem} \label{T0103}
Fix a positive real number $c_0>0$. 
For any positive real number $b>0$, there exists a positive integer 
$n_0$ such that $c_0\leq \log n_0(q-1)$, and for any integer 
$\disp{k=\frac{n(q-1)}{2q}(\log n(q-1)-c)}$ with $0\leq c\leq c_0$, we have 
\begin{equation} \label{0106}
\Vert \nu _n^{*k}-\pi _n\Vert _{TV}\geq 1-(4q+b)e^{-c}. 
\end{equation} 
\end{theorem}

The above theorem says that we can take a function $1-4qe^{-c}$ 
as a lower function of $\Vert \nu _n^{*k}-\pi _n\Vert _{TV}$. 

\section{Hamming schemes.}
\quad\,\,
In this section, we give notations of Hamming schemes, 
refering to \cite{BI}, \cite{CST}, \cite{D1}. 

For a positive integer $m$, we denote by $[m]=\{1, 2, \dots , m\}$, 
$[m]_0=\{0, 1, \dots , m-1\}$ and by $\sharp S$ the cardinal number 
for a finite set $S$. 

Let $n$ be a positive integer and $q$ an integer with $q\geq 2$. 
Put  
\begin{equation} \label{0201}
H(n, q)=[q]_0^n=\{x=(x_1, \dots , x_n)\,;\,x_j\in [q]_0\,(1\leq j\leq n)\}. 
\end{equation}
Take $x=(x_1, \dots , x_n)$, $y=(y_1, \dots , y_n)\in H(n, q)$ and 
define $d(x, y)$ by 
\begin{equation} \label{0202}
d(x, y)=\sharp \{j\in [n]\,;\,x_j\not= y_j\}. 
\end{equation}
For $x, x'\in H(n, q)$, we say that $x'$ is {\it{adjacent}} to $x$ 
if $d(x, x')=1$ and write $x\sim x'$. 
We fix an integer $q$ with $q\geq 2$ and denote by $X_n=H(n, q)$ 
for simplicity. 
We call an unordered pair $\{x, x'\}\subset X_n$ such that $x\sim x'$ 
the {\it{edge}} of $X_n$. 
Put 
\begin{equation} \label{0203}
E_{X_n}=\{\{x, x'\}\in X_n\,;\,x\sim x'\}. 
\end{equation}
Then the pair $(X_n, E_{X_n})$ is a simple undirected finite graph 
without loops. 
We call $H(n, q)$ the {\it{Hamming scheme}}. 
$X_n$ and $E_{X_n}$ are called the {\it{vertex set}} and 
the {\it{edge set}} of $H(n, q)$, respectively. 
We see that the distance $d(\cdot , \cdot )$ on $X_n$ coincides 
that derived from $E_{X_n}$. 

We denote by $S_m$ the symmetric group on $[m]$ or $[m]_0$ 
for a positive integer $m$. 
Let $G_n=S_q\wr S_n=S_q^n\rtimes S_n$ denote the {\it{wreath product}} 
of $S_q$ by $S_n$ with the product 
\begin{equation} \label{0204}
(\tau _1, \dots , \tau _n; \sigma )(\tau '_1, \dots , \tau '_n; \sigma ')
=(\tau _1\tau '_{\sigma ^{-1}(1)}, \dots , \tau _n\tau '_{\sigma ^{-1}(n)}; 
\sigma \sigma '), 
\end{equation}
where $\tau _j, \tau '_j\in S_q$\,($1\leq j\leq n)$, 
$\sigma , \sigma '\in S_n$, and we regard $S_q$ and $S_n$ as symmetric groups 
acting on $[q]_0$ and on $[n]$, respectively. 
$G_n$ acts on $X_n$ by 
\begin{equation} \label{0205}
g\cdot x=(\tau _1(x_{\sigma ^{-1}(1)}), \dots , \tau _n(x_{\sigma ^{-1}(n)})), 
\end{equation}
where $g=(\tau _1, \dots , \tau _n; \sigma )\in G_n$ and 
$x=(x_1, \dots , x_n)\in X_n$. 
The action of $G_n$ on $X_n$ is transitive. 
Put $x^{(0)}=(0, \dots , 0)\in X_n$. 
Then the stabilizer $H_n$ of $G_n$ at $x^{(0)}$ is given by 
\begin{equation} \label{0206}
H_n=S_{q-1}\wr S_n=\{(\tau _1, \dots , \tau _n; \sigma )\in G_n\,;\,
\tau _j\in S_{q-1}\,\,{\text{for all}}\,\,j\in [n]\}, 
\end{equation}
where we regard $S_{q-1}$ as a symmetric group acting on $[q-1]$. 
For each integer $j\in \{0, 1, \dots , n\}$ we put 
\begin{equation} \label{0207}
x^{(j)}=(\overbrace{1, \dots , 1}^j, 0, \dots , 0)\in X_n. 
\end{equation}
Then, $X_n=G_n/H_n$ and we have the $H_n$-orbit decomposion 
\begin{equation} \label{0208}
X_n=\bigcup _{j=0}^nH_n\cdot x^{(j)}. 
\end{equation}
We see that $x\sim x'$ implies that $g\cdot x\sim g\cdot x'$ 
for any $x, x'\in X_n$ and $g\in G_n$. 
Hence, for any $j\in \{0, 1, \dots , n\}$, we have that 
$x\in H_n\cdot x^{(j)}$ if and only if $d(x^{(0)}, x)=j$. 
For $j\in \{0, 1, \dots , n\}$, put 
\begin{equation} \label{0209}
g^{(j)}=(\overbrace{(0, 1), \dots , (0, 1)}^j, 1_{S_q}, \dots , 1_{S_q}; 
1_{S_n}), 
\end{equation}
where $(0, 1)\in S_q$ is the transposition of $0$ and $1$, and 
$1_{S_q}\in S_q$, $1_{S_n}\in S_n$ are the identity permutations. 
We see that $g^{(j)}\cdot x^{(0)}=x^{(j)}$. 
Hence, we have the decomposition $\disp{G_n=\bigcup _{j=0}^n
H_ng^{(j)}H_n}$ of $G_n$ into $H_n$-double cosets. 
Let $L^1(G_n)$ denote the algebra of all functions on $G_n$ 
with the convolution 
\begin{equation} \label{0210}
f_1*f_2(g)=\sum _{g'\in G_n}f_1(g(g')^{-1})f_2(g')
=\sum _{g'\in G_n}f_1(g')f_2((g')^{-1}g), 
\end{equation}
where $f_1, f_2\in L^1(G_n)$ and $g\in G_n$. 
We see that $(G_n, H_n)$ is a {\it{Gelfand pair}}, that is, 
the subalgebra $L^1(H_n\backslash G_n/H_n)\subset L^1(G_n)$ 
of all $H_n$-biinvariant functions on $G_n$ is a commutative algebra 
since $(g^{(j)})^{-1}=g^{(j)}$ for any $j\in \{0, 1, \dots , n\}$\,
(see \cite{CST}, Example 4.3.2). 
We denote by $L(X_n)$ the Hilbert space of all functions on $X_n$ 
with the inner product 
\begin{equation} \label{0211}
\langle f_1, f_2\rangle _{L(X_n)}=\sum _{x\in X_n}f_1(x)\overline{f_2(x)}, 
\end{equation}
where $f_1, f_2\in L(X_n)$, and write $\Vert f \Vert _{L(X_n)}
=\langle f, f\rangle _{L(X_n)}^{\frac{1}{2}}$ for $f\in L(X)$. 
$G_n$ acts on $L(X_n)$ by 
\begin{equation} \label{0212}
(g\cdot f)(x)=f(g^{-1}\cdot x), 
\end{equation}
where $g\in G_n$, $f\in L(X_n)$ and $x\in X_n$. 
It is easy to show that the action is unitary. 

Let $W$ be a $G_n$-module. 
We denote by $W_{H_n}$ the subspace of all $H_n$-invariant elements in $W$, 
that is, 
\begin{equation} \label{0213}
W_{H_n}=\{w\in W\,;\,h\cdot w=w\,\,{\text{for all}}\,\,h\in H_n\}. 
\end{equation}

The condition that $(G_n, H_n)$ is a Gelfand pair indicates the properties 
of irreducible components appearing in $L(X_n)$. 

\begin{lemma} \label{T0201}
Let $\disp{L(X_n)=\bigoplus _{\lambda \in \Lambda}V_{\lambda}}$ be 
an irreducible decomposition of $L(X_n)$.
\newline
{\rm{(1)}}\,\,
$L(X_n)$ is multiplicity-free, that is, $V_{\lambda}$ is not equivalent 
to $V_{\lambda '}$ if $\lambda \not= \lambda '$. 
\newline
{\rm{(2)}}\,\,
For any $\lambda \in \Lambda$, we have 
$\dim (V_{\lambda})_{H_n}=1$. 
\end{lemma}

{\it{Proof}}.\,\,
See \cite{CST}, Theorem 4.4.2, 4.6.2 for example. 
\qed

\vspace{1pc}

Each irreducible component $V_{\lambda}$ is called 
the {\it{spherical representation}} for $(G_n, H_n)$. 

We construct irreducible components in $L(X_n)$. 
Take an integer $a\in [q]_0$ and define a function 
$\chi _a\,:\,[q]_0\longrightarrow {\mathbb{C}}$ by 
\begin{equation*}
\chi _a(x)=\zeta _q^{ax}, 
\end{equation*}
where $x\in [q]_0$ and $\disp{\zeta _q=\exp \frac{2\pi i}{q}}\in {\mathbb{C}}$ 
is a primitive $q$-th root of $1$ in ${\mathbb{C}}$.
For $x\in [q]_0$, we have 
\begin{equation*}
\sum _{a=0}^{q-1}\chi _a(x)=\begin{cases}
q, & x=0,      \\
0, & x\not= 0. 
\end{cases}
\end{equation*}
For $a=(a_1, \dots , a_n)\in X_n$, we define a function 
$\chi _a\,:\,X_n\longrightarrow {\mathbb{C}}$ by 
\begin{equation} \label{0214}
\chi _a(x)=\prod _{j=1}^n\chi _{a_j}(x_j)=\zeta _q^{a_1x_1+\cdots +a_nx_n}, 
\end{equation}
where $x=(x_1, \dots , x_n)\in X_n$. 
Then, we see that $\{\chi _a\in L(X_n)\,;\,a\in X_n\}$ is 
an orthogonal basis for $L(X_n)$ and 
$\Vert \chi _a\Vert _{L(X_n)}=q^{\frac{n}{2}}$ for $a\in X_n$. 

For $j\in \{0, 1, \dots , n\}$, we put 
\begin{equation} \label{0215}
V_j=\bigoplus _{\sharp \{l\in [n]\,;\,a_l\not= 0\}=j}{\mathbb{C}}\chi _a. 
\end{equation}
Then, $V_j$ is $G_n$-invariant for each $j$ with $0\leq j\leq n$ and 
we have the orthogonal desomposition 
\begin{equation} \label{0216}
L(X_n)=\bigoplus _{j=0}^nV_j. 
\end{equation}

For $j\in \{0, 1, \dots , n\}$, we put 
\begin{equation} \label{0217}
\omega _j=\sum _{\sharp \{l\in [n]\,;\,a_l\not= 0\}=j}\chi _a. 
\end{equation}
Then, $\omega _j$ is nonzero $H_n$-invariant element in $V_j$, and 
any $H_n$-invariant element in $V_j$ is the scalar multiple of $\omega _j$. 
Hence, all $V_j$'s are irreducible and $L(X_n)$ is multiplicity-free. 

For $j\in \{0, 1, \dots , n\}$, we define a function 
$\phi _j\,:\,G_n\longrightarrow {\mathbb{C}}$ by 
\begin{equation} \label{0218}
\phi _j(g)=\left\langle \frac{\omega _j}{\Vert \omega _j\Vert_{L(X_n)}},\,
g\cdot \frac{\omega _j}{\Vert \omega _j\Vert_{L(X_n)}}\right\rangle _{L(X_n)}
=\frac{1}{\Vert \omega _j\Vert _{L(X_n)}^2}
\langle \omega _j, g\cdot \omega _j\rangle _{L(X_n)}. 
\end{equation}
Then, $\phi _j$ is $H_n$-biinvariant and $\phi _j(1_{G_n})=1$, 
where $1_{G_n}\in G_n$ is the unit element. 
Moreover, $\phi _j$ is real-valued since $(g^{(j)})^{-1}\in H_ng^{(j)}H_n$ 
for any $j\in \{0, 1, \dots , n\}$\,(see \cite{CST}, Theorem 4.8.2). 
$\phi _j$ is called the {\it{spherical function}} on $G_n$. 
We regard $\phi _j$ as an $H_n$-invariant function on $X_n$.
$\phi _j$ is calculated as 
\begin{equation} \label{0219}
\phi _j(g^{(l)})=\frac{1}{\disp{\binom{n}{j}}}\sum _{r=0}^j\binom{l}{r}
\binom{n-l}{j-r}\left( -\frac{1}{q-1}\right) ^r, 
\end{equation}
where $l\in \{0, 1, \dots , n\}$\,
(see \cite{CST}, Theorem 5.3.2). 
We give another realization of $\phi _j(g^{(l)})$. 
For a complex number $\alpha$ and a nonnegative integer $m$, put
\begin{equation*}
(\alpha )_m=\begin{cases}
\alpha (\alpha +1)\cdots (\alpha +m-1), & m\geq 1, \\
1,                                      & m=0. 
\end{cases}
\end{equation*}
$(\alpha )_m$ is called the {\it{Pochhammer symbol}}. 
Take complex numbers $\alpha , \beta , \gamma \in {\mathbb{C}}$, a variable $x$, 
and define the Gauss hypergeometric series 
\begin{equation*}
F\left( \begin{array}{c}
\alpha , \beta \\
\gamma
\end{array}; x\right) =\sum _{m=0}^{\infty}\frac{(\alpha )_m(\beta )_m}
{(\gamma )_mm!}x^m. 
\end{equation*}
We write $\phi _j(l)=\phi _j(g^{(l)})$ for simplicity. 
The polynomial $\phi _j$ is called the {\it{Krawtchouk polynomial}}. 
Using a Gauss hypergeometric series, $\phi _j(l)$ is realized as 
\begin{equation} \label{0220}
\phi _j(l)=F\left( \begin{array}{c}
-j, -l \\
-n
\end{array}; \frac{q}{q-1}\right) 
=\sum _{r=0}^j\frac{(-j)_r(-l)_r}{(-n)_rr!}\left( \frac{q}{q-1}\right) ^r, 
\end{equation}
where $j, l\in \{0, 1, \dots , n\}$. 

For $f\in L^1(H_n\backslash G_n/H_n)$, we define 
\begin{equation} \label{0221}
\widehat{f}(\phi _j)=\sum _{g\in G_n}f(g)\overline{\phi _j(g)}
=\sum _{g\in G_n}f(g)\phi _j(g), 
\end{equation}
where $j\in \{0, 1, \dots , n\}$. 
$\widehat{f}$ is called the {\it{spherical transform}} 
of $f\in L^1(H_n\backslash G_n/H_n)$. 
For $f_1, f_2\in L^1(H_n\backslash G_n/H_n)$, We have 
\begin{equation} \label{0222}
(f_1*f_2)\,\widehat{}=\widehat{f}_1\widehat{f}_2. 
\end{equation}

Take $f\in L(X_n)_{H_n}$, $j\in \{0, 1, \dots , n\}$ and define 
\begin{equation} \label{0223}
{\cal{F}}(f)(\phi _j)=\sum _{x\in X_n}f(x)\overline{\phi _j(x)}
=\sum _{x\in X_n}f(x)\phi _j(x). 
\end{equation}
${\cal{F}}(f)$ is called the {\it{spherical transform}} of $f\in L(X_n)_{H_n}$. 

For an $H_n$-invariant function $f\in L(X_n)_{H_n}$ on $X_n$, we denote 
by $\widetilde{f}$ the $H_n$-biinvariant function on $G_n$ corresponding to $f$. 
We see that 
\begin{align*}
{\cal{F}}(f)(\phi _j) & 
=\sum _{x\in X_n}f(x)\phi _j(x)
=\sum _{x\in X_n}\frac{1}{\sharp H_n}\sum _{g\cdot x^{(0)}=x}
f(g\cdot x^{(0)})\phi _j(g) \\
{} & 
=\frac{1}{\sharp H_n}\sum _{g\in G_n}\widetilde{f}(g)\phi _j(g)
=\frac{1}{\sharp H_n}(\widetilde{f})\,\widehat{}\,(\phi _j). 
\end{align*} 
Take $f_1, f_2\in L(X_n)_{H_n}$ and define $f_1*f_2\in L(X_n)_{H_n}$ such that 
\begin{equation} \label{0224}
(f_1*f_2)\,\widetilde{}=\frac{1}{\sharp H_n}(\widetilde{f}_1*\widetilde{f}_2). 
\end{equation}

\begin{lemma} \label{T0202}
Let $f_1, f_2\in L(X_n)_{H_n}$ be two elements in $L(X_n)_{H_n}$. 
For $j\in \{0, 1, \dots , n\}$, we have 
\begin{equation} \label{0225}
{\cal{F}}(f_1*f_2)={\cal{F}}(f_1){\cal{F}}(f_2). 
\end{equation}
\end{lemma}

{\it{Proof}}.\,\,
By (\ref{0222}) and (\ref{0223}), we have 
\begin{align*}
{\cal{F}}(f_1*f_2)(\phi _j) & 
=\frac{1}{\sharp H_n}((f_1*f_2)\,\widetilde{}\,)\,\widehat{}\,(\phi _j)
=\frac{1}{(\sharp H_n)^2}((\widetilde{f}_1)*(\widetilde{f}_2))\,
\widehat{}\,(\phi _j) \\
{} & 
=\frac{1}{(\sharp H_n)^2}((\widetilde{f}_1)\,\widehat{}\,
(\widetilde{f}_2)\,\widehat\,)(\phi _j)
={\cal{F}}(f_1)(\phi _j){\cal{F}}(f_2)(\phi _j). 
\end{align*}
\qed

\vspace{1pc}

For $f\in L(X_n)_{H_n}$ and a nonnegative integer $k$, we define 
$f^{*k}$ recursively by 
\begin{equation} \label{0226}
f^{*0}=\delta _{x^{(0)}}, \quad f^{*k}=f^{*(k-1)}*f,\,\,k\geq 1. 
\end{equation}

\begin{lemma} \label{T0203}
For $j\in \{0, 1, \dots , n\}$ and a nonnegative integer $k$, we have 
\begin{equation} \label{0227}
{\cal{F}}(f^{*k})(\phi _j)={\cal{F}}(f)(\phi _j)^k. 
\end{equation}
\end{lemma}

{\it{Proof}}.\,\,
We prove it by induction in $k$. 
We see that 
\begin{equation*}
{\cal{F}}(f^{*0})(\phi _j)=\sum _{x\in X_n}\delta _{x^{(0)}}(x)\phi _j(x)
=\phi _j(x^{(0)})=1. 
\end{equation*}
We assume $k\geq 1$ and the claim satisfies for any integer less than $k$. 
Then we have 
\begin{equation*}
{\cal{F}}(f^{*k})(\phi _j)={\cal{F}}(f^{*(k-1)}*f)
={\cal{F}}(f^{*(k-1)}){\cal{F}}(f)={\cal{F}}(f)^{k-1}{\cal{F}}(f)
={\cal{F}}(f)^k.  
\end{equation*}
\qed

\vspace{1pc}

Take $x, x'\in X_n$ and define the {\it{transition probability}} $p_n(x, x')$ by 
\begin{equation} \label{0228}
p_n(x, x')=\begin{cases}
\disp{\frac{1}{n(q-1)}}, & x\sim x',     \\
\quad\,\,\,\,0,                       & x\not\sim x'. 
\end{cases}
\end{equation}
Since $G_n$ preserves adjacency on $X_n$, for $x, x'\in X_n$ and $g\in G_n$, 
we see that 
\begin{equation} \label{0229}
p_n(g\cdot x, g\cdot x')=p_n(x, x'). 
\end{equation}
Put 
\begin{equation} \label{0230}
\nu _n(x)=p_n(x^{(0)}, x). 
\end{equation}
Then, $\nu _n$ is an $H_n$-invariant probability measure on $X_n$. 

\begin{lemma} \label{T0204}
For $j\in \{0, 1, \dots , n\}$, we have 
\begin{align} \label{0231}
{\cal{F}}(\nu _n)(\phi _j)=1-\frac{jq}{n(q-1)}.
\end{align}
\end{lemma}

{\it{Proof}}.\,\,
Since $x^{(0)}\sim x$ if and only if $x\in H_n\cdot x^{(1)}$, we see that 
\begin{equation*}
\phi _j(x^{(1)})=F\left( \begin{array}{c}
-j, -1 \\
-n
\end{array}; \frac{q}{q-1}\right) 
=1+\frac{(-j)\cdot (-1)}{(-n)\cdot 1}\cdot \frac{q}{q-1}
=1-\frac{jq}{n(q-1)}. 
\end{equation*}
Hence, we have 
\begin{align*}
{\cal{F}}(\nu _n)(\phi _j) & 
=\sum _{x\in X_n}\nu _n(x)\phi _j(x)
=\sum _{x^{(0)}\sim x}\nu _n(x)\phi _j(x)
=\sum _{x\in H_n\cdot x^{(1)}}\nu _n(x)\phi _j(x^{(1)}) \\
{} & 
=n(q-1)\cdot \frac{1}{n(q-1)}\left( 1-\frac{jq}{n(q-1)}\right) 
=1-\frac{jq}{n(q-1)}. 
\end{align*}
\qed

\section{Upper bounds.}
\quad\,\,
In this section, we give a majorant function for the distances 
between the $k$-step transitions distributions and the distributions 
of equilibrium for the simple random walks on the Hamming schemes 
$\{H(n, q)\}$ with $q\geq 3$.  

We denote by $M(X_n)$ the vector space of all complex-valued measure on $X_n$. 
For a measure $\mu \in M(X_n)$ on $X_n$, we put 
\begin{equation} \label{0301}
\Vert \mu \Vert _{TV}=\max \{|\mu (S)|\,;\,S\subset X_n\}. 
\end{equation}
For two measure $\mu , \nu \in M(X_n)$, we define the {\it{total variance distance}}
by $\Vert \mu -\nu \Vert _{TV}$. 
We regard a measure $\mu \in M(X_n)$ on $X_n$ as a function 
$\mu\,:\,X_n\longrightarrow {\mathbb{C}}$ defined by 
$\mu (x)=\mu (\{x\})$ for $x\in X_n$. 
Take two probability measures $\mu , \nu \in M(X_n)$ on $X_n$. 
Then we have an equality 
\begin{equation*}
\Vert \mu -\nu \Vert _{TV}^2=\left( \frac{1}{2}\sum _{x\in X_n}
|\mu (x)-\nu (x)|\right) ^2
\leq \frac{\sharp X_n}{4}\Vert \mu -\nu \Vert _{L(X_n)}^2. 
\end{equation*}
We denote by $\pi _n$ the uniform probability measure on $X_n$, that is, 
\begin{equation} \label{0302}
\pi _n(S)=\frac{\sharp S}{\sharp X_n}=\frac{\sharp S}{q^n}, 
\end{equation}
where $S\subset X_n$ is a subset of $X_n$. 
We find the upper bound of the total variance distance 
$\Vert \nu _n^{*k}-\mu _n\Vert _{TV}$ 
with the Fourier transforms of spherical functions. 

\begin{lemma}[The upper bound lemma] \label{T0301}
For a nonnegative integer $k$, we have 
\begin{equation} \label{0303}
\Vert \nu _n^{*k}-\pi _n\Vert _{TV}^2\leq \frac{1}{4}
\sum _{j=1}^nd_j|{\cal{F}}(\nu _n)(\phi _j)|^2. 
\end{equation}
\end{lemma}

{\it{Proof}}.\,\,
See \cite{D1},  Chapter 3B, Lemma 1 or \cite{CST}, Corollary 4.9.2. 
\qed

\vspace{1pc}

Before estimating the total variance distance 
$\Vert \nu _n^k-\pi _n\Vert _{TV}$, we give two inequalities. 

\begin{lemma} \label{T0302}
{\rm{(1)}}\,\,
For a real number $x$ such that $\disp{x\leq \frac{5}{4}}$, 
we have $e^{-x}\geq |1-x|$. 
\newline
{\rm{(2)}}\,\,
If the real number $x$ satisfies the condition $\disp{x\geq \frac{4}{3}}$, 
one has $e^{-x}\leq |1-x|$. 
\end{lemma}

{\it{Proof}}.\,\,
First, we see that 
\begin{equation*}
\frac{8}{3}=\sum _{j=0}^3\frac{1}{j!}\leq e=\sum _{j=0}^{\infty}\frac{1}{j!}
\leq \frac{5}{2}+\sum _{j=3}^{\infty}\frac{1}{2\cdot 3^{j-2}}
=\frac{5}{2}+\frac{1}{4}=\frac{11}{4}. 
\end{equation*}
(1)\,\,
By Taylor's theorem, for $x\in {\mathbb{R}}$, there exists a real number 
$\theta \in {\mathbb{R}}$ with $0<\theta <1$ such that 
\begin{equation*}
e^{-x}=1-x+\frac{x^2}{2}e^{-\theta x}. 
\end{equation*}
If $x\leq 1$, we have $e^{-x}\geq 1-x=|1-x|$. 
Hence, we only consider the case $\disp{1\leq x\leq \frac{5}{4}}$. 
We see that 
\begin{equation*}
e^5\leq \left( \frac{11}{4}\right) ^5=\frac{161051}{1024}
\leq 256=4^4. 
\end{equation*}
So $\disp{e^{-\frac{5}{4}}\geq \frac{1}{4}=\frac{5}{4}-1}$. 
Therefore, for any real number $x$ such that $\disp{1\leq x\leq \frac{5}{4}}$, 
we have 
\begin{equation*}
e^{-x}\geq e^{-\frac{5}{4}}\geq \frac{5}{4}-1\geq x-1=|1-x|. 
\end{equation*}
(2)\,\,
We see that 
\begin{equation*}
e^4\geq \left( \frac{8}{3}\right) ^4=\frac{4096}{81}\geq 27=3^3. 
\end{equation*}
Hence, $\disp{e^{-\frac{4}{3}}\leq \frac{1}{3}=\frac{4}{3}-1}$. 
This implies that for any real number $x$ with $\disp{x\geq \frac{4}{3}}$, 
\begin{equation*}
e^{-x}\leq e^{-\frac{4}{3}}\leq \frac{4}{3}-1\leq x-1=|1-x|. 
\end{equation*}
\qed

\vspace{1pc}

Here, we give a majorant function for the total variance distance 
$\Vert \nu _n^{k}-\pi _n\Vert _{TV}$ with a large integer $q$. 

\begin{theorem} \label{T0303}
Let $q$ be an integer with $q\geq 5$. 
Take a positive integer $k$ such that $\disp{k=\frac{n(q-1)}{2q}
(\log n(q-1)+c)}$ with $c>0$. 
Then, we have 
\begin{equation*}
\Vert \nu _n^{*k}-\pi _n\Vert _{TV}^2
\leq \frac{1}{4}(e^{e^{-c}}-1). 
\end{equation*}
\end{theorem}

{\it{Proof}}.\,\,
The condition $q\geq 5$ implies that $\disp{\frac{jq}{n(q-1)}
\leq \frac{q}{q-1}\leq \frac{5}{4}}$ for any integer 
$j\in \{0, 1, \dots , n\}$. 
Hence, by Lemma \ref{T0301}, we have 
\begin{align*}
\Vert \nu _n^{*k}-\pi _n\Vert _{TV}^2 & 
\leq \frac{1}{4}\sum _{j=1}^nd_j|{\cal{F}}(\nu _n)(\phi _j)|^{2k}
=\frac{1}{4}\sum _{j=1}^n\binom{n}{j}(q-1)^j
\left| 1-\frac{jq}{n(q-1)}\right| ^{2k} \\
{} & 
\leq \frac{1}{4}\sum _{j=1}^n\frac{n^j(q-1)^j}{j!}e^{-\frac{2jkq}{n(q-1)}}
=\frac{1}{4}\sum _{j=1}^n\frac{1}{j!}e^{j\left( \log n(q-1)
-\frac{2kq}{n(q-1)}\right)}. 
\end{align*}
Since $\disp{k=\frac{n(q-1)}{2q}(\log n(q-1)+c)}$, we have 
\begin{equation*}
\Vert \nu _n^{*k}-\pi _n\Vert _{TV}^2
\leq \frac{1}{4}\sum _{j=1}^n\frac{e^{-cj}}{j!}
\leq \frac{1}{4}\sum _{j=1}^{\infty}\frac{e^{-cj}}{j!}
=\frac{1}{4}(e^{e^{-c}}-1). 
\end{equation*}
\qed

\vspace{1pc}

It remains to show the case $q=3, 4$. 
For any real number $\alpha \in {\mathbb{R}}$, we denote by 
\begin{align*}
\lfloor \alpha \rfloor & 
=\max \{m\in {\mathbb{Z}}\,;\,m\leq \alpha \}, \\
\lceil \alpha \rceil & 
=\min \{m\in {\mathbb{Z}}\,;\,m\geq \alpha \}. 
\end{align*} 
In order to estimate a total variance distances, where $q=3, 4$, we give a lemma. 

\begin{lemma} \label{T0304}
{\rm{(1)}}\,\,
Let $m$ be an integer with $m\geq 2$ and $l$ an integer 
such that $0\leq l\leq m$. 
Put $\disp{f_m(x)=\frac{x-m+2}{x+m}=1-\frac{2m-2}{x+m}}$. 
Then, we have 
\begin{equation}
\sum _{p=l}^{2m-l-1}f_m(p)\leq 2\log \frac{3m-l}{l+m}\leq \log 9.  
\end{equation}
{\rm{(2)}}\,\,
Assume $m$ is an integer with $m\geq 2$ and $l$ is an integer such that 
$\disp{0\leq l\leq \left\lfloor \frac{m}{2}\right\rfloor}$. 
Put $\disp{f_m(x)=\frac{2x-m+3}{x+m}=2-\frac{3m-3}{x+m}}$. 
Then, one has 
\begin{equation}
\sum _{p=l}^{m-l-1}f_m(p)\leq 3\log \frac{2m-l}{l+m}\leq \log 8. 
\end{equation}
\end{lemma}

{\it{Proof}}.\,\,
(1)\,\,
$f_m$ is a monotone incresing continuous function 
on the open interval $(-m, +\infty)$. 
Hence, 
\begin{align*}
\sum _{p=l}^{2m-l-1}f_m(p) & 
\leq \int _l^{2m-l}f_m(x)dx
=\int _l^{2m-l}\left( 1-\frac{2m-2}{x+m}\right) dx \\
{} & 
=2m-2l-(2m-2)\log \frac{3m-l}{l+m}. 
\end{align*}
Put 
\begin{equation*}
g_m(x)=x+m\log \frac{3m-x}{x+m}. 
\end{equation*}
If $0\leq x\leq m$, we have 
\begin{equation*}
g'_m(x)=1+m\left( \frac{-1}{3m-x}-\frac{1}{x+m}\right) 
=\frac{-(x-m)^2}{(3m-x)(x+m)}\leq 0. 
\end{equation*}
Hence, $g_m(l)\geq g_m(m)=m$ for $l\in \{0, 1, \dots , m\}$. 
Therefore, 
\begin{equation*}
\sum _{p=l}^{m-l-1}f_m(p)\leq 2m-2g_m(l)+2\log \frac{3m-l}{l+m}
\leq 2\log \frac{3m-l}{l+m}\leq 2\log 3=\log 9. 
\end{equation*}
(2)\,\,
Similarly to (1), we have 
\begin{align*}
\sum _{p=l}^{m-l-1}f_m(p) & 
\leq \int _l^{m-l}f_m(x)dx
=\int _l^{m-l}\left( 2-\frac{3m-3}{x+m}\right) dx \\
{} & 
=2m-4l-(3m-3)\log \frac{2m-l}{l+m}. 
\end{align*}
Put 
\begin{equation*}
g_m(x)=4x+3m\log \frac{2m-x}{x+m}. 
\end{equation*}
For a real number $x$ with $\disp{0\leq x\leq \frac{m}{2}}$, we see that  
\begin{align*}
g'_m(x) & 
=4+3m\left( \frac{-1}{2m-x}-\frac{1}{x+m}\right) 
=\frac{-(2x-m)^2}{(2m-x)(x+m)}\leq 0. 
\end{align*}
Hence, $\disp{g_m(l)\geq g_m\left( \frac{m}{2}\right) =2m}$ 
for any integer $l$ such that $\disp{0\leq l\leq \left\lfloor 
\frac{m}{2}\right\rfloor}$. 
Therefore, 
\begin{equation*}
\sum _{p=l}^{m-l-1}f_m(p)\leq 3\log \frac{2m-l}{l+m}+2m-g_m(l)
\leq 3\log \frac{2m-l}{l+m}\leq 3\log 2=\log 8. 
\end{equation*}
\qed

\vspace{1pc}

\begin{lemma}\label{T0305}
{\rm{(1)}}\,\,
Put $\disp{a_{n, j}=2^j\binom{n}{j}}$ for integers $n$ and $j$ 
such that $n\geq 3$ and that $0\leq j\leq n$. 
Then, for any integers $m$ and $l$ such that $m\geq 2$ and that 
$0\leq l\leq m-1$, we have 
\begin{equation} \label{0304}
\frac{a_{3m-3, 3m-l-3}}{a_{3m-3, l+m-1}}
\leq \frac{a_{3m-2, 3m-l-2}}{a_{3m-2, l+m-1}}
\leq \frac{a_{3m-1, 3m-l-1}}{a_{3m-1, l+m-1}}
\leq 9. 
\end{equation}
{\rm{(2)}}\,\,
We put $\disp{a_{n, j}=3^j\binom{n}{j}}$ for integers $n$ and $j$ 
such that $n\geq 2$ and that $0\leq j\leq n$. 
In this case, for any integer $m$ and $l$ such that $m\geq 2$ and 
that $\disp{0\leq l\leq \left\lfloor \frac{m-1}{2}\right\rfloor}$, we see that 
\begin{equation} \label{0305}
\frac{a_{2m-2, 2m-l-2}}{a_{2m-2, l+m-1}}
\leq \frac{a_{2m-1, 2m-l-1}}{a_{2m-1, l+m-1}}
\leq 8. 
\end{equation}
\end{lemma}

{\it{Proof}}.\,\,
(1)\,\,
Assume $n=3m-1$, where $m$ is an integer such that $m\geq 2$. 
For an integer $l$ such that $0\leq l\leq m-1$, we have 
\begin{align*}
\frac{a_{3m-1, 3m-l-1}}{a_{3m-1, l+m-1}} & 
=2^{2m-2l}\frac{\disp{\binom{3m-1}{3m-l-1}}}{\disp{\binom{3m-1}{l+m-1}}}
=2^{2m-2l}\frac{(2m-l)(2m-l-1)\cdots (l+1)}{(3m-l-1)(3m-l-2)\cdots (l+m)} \\
{} & 
=\prod _{p=l}^{2m-l-1}\frac{2p+2}{p+m}
=\prod _{p=l}^{2m-l-1}\left(1+\frac{p-m+2}{p+m}\right) \\
{} & 
\leq \prod _{p=l}^{2m-l-1}\exp \left( \frac{p-m+2}{p+m}\right) 
=\exp \sum _{p=l}^{2m-l-1}\frac{p-m+2}{p+m}. 
\end{align*}
By Lemma \ref{T0304}\,(1), we have 
\begin{equation*}
\frac{a_{3m-1, 3m-l-1}}{a_{3m-1, l+m-1}}
\leq \exp (\log 9)=9. 
\end{equation*}
Next, we consider the case $n=3m-2$. 
We see that $2(2m-l)-(3m-l-1)=m-l+1>0$ for an integer 
$l$ with $0\leq l\leq m-1$. 
Hence, 
\begin{align*}
\frac{a_{3m-2, 3m-l-2}}{a_{3m-2, l+m-1}} & 
=2^{2m-2l-1}\frac{\disp{\binom{3m-2}{3m-l-2}}}{\disp{\binom{3m-2}{l+m-1}}}
=\frac{3m-l-1}{2(2m-l)}\cdot \frac{a_{3m-1, 3m-l-1}}{a_{3m-1, l+m-1}} \\
{} & 
\leq \frac{a_{3m-1, 3m-l-1}}{a_{3m-1, l+m-1}}. 
\end{align*}
Finally, we assume $n=3m-3$ for an integer $m$ such that $m\geq 2$. 
For $l\in \{0, 1, \dots , m-1\}$, we see that 
$2(2m-l-1)-(3m-l-2)=m-l>0$ and that 
\begin{align*}
\frac{a_{3m-3, 3m-l-3}}{a_{3m-3, l+m-1}} & 
=2^{2m-2l-2}\frac{\disp{\binom{3m-3}{3m-l-3}}}{\disp{\binom{3m-3}{l+m-1}}}
=\frac{3m-l-2}{2(2m-l-1)}\frac{a_{3m-1, 3m-l-2}}{a_{3m-2, l+m-1}} \\ 
{} & 
\leq \frac{a_{3m-2, 3m-l-2}}{a_{3m-2, l+m-1}}.
\end{align*}
(2)\,\,
Assume $n=2m-1$ for some integer $m$ with $m\geq 2$. 
For an integer $j$ such that $\disp{0\leq l\leq \left\lfloor 
\frac{m-1}{2}\right\rfloor}$, we have 
\begin{align*}
\frac{a_{2m-1, 2m-l-1}}{a_{2m-1, l+m-1}} & 
=3^{m-2l}\frac{\disp{\binom{2m-1}{2m-l-1}}}{\disp{\binom{2m-1}{l+m-1}}}
=3^{m-2l}\frac{(m-l)(m-l-1)\cdots (l+1)}{(2m-l-1)(2m-l-2)\cdots (l+m)} \\
{} & 
=\prod _{p=l}^{m-l-1}\frac{3p+3}{p+m}
=\prod _{p=l}^{m-l-1}\left( 1+\frac{2p-m+3}{p+m}\right) \\
{} & 
\leq \prod _{p=l}^{m-l-1}\exp \left( \frac{2p-m+3}{p+m}\right) 
=\exp \sum _{p=l}^{m-l-1}\frac{2p-m+3}{p+m}. 
\end{align*}
By Lemma \ref{T0304}\,(2), we have 
\begin{equation*}
\frac{a_{2m-1, 2m-l-1}}{a_{2m-1, l+m-1}}\leq \exp (\log 8)=8. 
\end{equation*}
Next, we consider the case $n=2m-2$ for some integer $m$ with $m\geq 2$. 
For an integer $l$ such that $\disp{0\leq l\leq \left\lfloor 
\frac{m-1}{2}\right\rfloor}$, we have 
\begin{align*}
\frac{a_{2m-2, 2m-l-2}}{a_{2m-2, l+m-1}} & 
=3^{m-2l-1}\frac{\disp{\binom{2m-2}{2m-l-2}}}{\disp{\binom{2m-2}{l+m-1}}}
=\frac{2m-l-1}{3(m-l)}\cdot \frac{a_{2m-1, 2m-l-1}}{a_{2m-1, l+m-1}} \\
{} & 
\leq \frac{a_{2m-1, 2m-l-1}}{a_{2m-1, l+m-1}}\leq 8, 
\end{align*}
since $3(m-l)-(2m-l-1)=m-2l+1\geq 0$. 
\qed

\vspace{1pc}

Using these lemmas, we estimate total variance distance 
$\Vert \nu _n^{*k}-\pi _n\Vert _{TV}$ for $q=3, 4$. 

\begin{theorem} \label{T0306}
Assume $q=3$. 
For a positive integer $\disp{k=\frac{n}{3}(\log 2n+c)}$ 
with $n\geq 3$ and $c>0$, we have 
\begin{equation} \label{0306}
\Vert \nu _n^{*k}-\pi _n\Vert _{TV}^2
\leq \frac{5}{2}(e^{e^{-c}}-1). 
\end{equation}
\end{theorem}

{\it{Proof}}.\,\,
By Lemma \ref{T0301}, we see that 
\begin{equation*}
\Vert \nu _n^{*k}-\pi _n\Vert _{TV}^2
\leq \frac{1}{4}\sum _{j=1}^nd_j|{\cal{F}}(\nu _n)(\phi _j)|^{2k}
=\frac{1}{4}\sum _{j=1}^n2^j\binom{n}{j}\left| 1-\frac{3j}{2n}\right| ^{2k}. 
\end{equation*}
Put $\disp{a_{n, j}=2^j\binom{n}{j}}$. 
Assume $n=3m-1$, where $m$ is an integer such that $m\geq 2$. 
By Lemma \ref{T0305}\,(1), for an integer $l$ such that $0\leq l\leq m-1$, 
we have $\disp{\frac{a_{3m-1, 3m-l-1}}{a_{3m-1, l+m-1}}\leq 9}$. 
Moreover, for an integer $l$ such that $0\leq l\leq m-1$, we have 
$3(3m-l-1)-2(3m-1)=3m-3l-1>0$, $2(3m-1)-3(l+m-1)=3m-3l+1>0$ and 
\begin{align*}
\left| 1-\frac{3(3m-l-1)}{2(3m-1)}\right| & 
=\frac{3(3m-l-1)}{2(3m-1)}-1=\frac{3m-3l-1}{2(3m-1)} \\
{} & 
\leq \frac{3m-3l+1}{2(3m-1)}=1-\frac{3(l+m-1)}{2(3m-1)}
=\left| 1-\frac{3(l+m-1)}{2(3m-1)}\right| . 
\end{align*}
Hence, 
\begin{align*}
\sum _{j=1}^n2^j\binom{n}{j}\left| 1-\frac{3j}{2n}\right| ^{2k} & 
\leq \sum _{j=1}^{m-2}a_{n, j}\left| 1-\frac{3j}{2n}\right| ^{2k}
+a_{n, 2m-1}\left| 1-\frac{3(2m-1)}{2n}\right| ^{2k} \\
{} & \quad 
+\sum _{l=0}^{m-1}(a_{3m-1, l+m-1}+a_{3m-1, 3m-l-1})
\left| 1-\frac{3(l+m-1)}{2(3m-1)}\right| ^{2k} \\
{} & 
\leq 10\sum _{j=1}^{2m-1}a_{n, j}\left| 1-\frac{3j}{2n}\right| ^{2k}
\leq 10\sum _{j=1}^{2m-1}\frac{2^jn^j}{j!}e^{-\frac{3jk}{n}} \\
{} & 
\leq 10\sum _{j=1}^{\infty}\frac{1}{j!}
e^{j\left( \log 2n-\frac{3k}{n}\right)}. 
\end{align*}
Next, we consider the case $n=3m-2$. 
We see that $\disp{\frac{a_{3m-2, 3m-l-2}}{a_{3m-2, l+m-1}}\leq 9}$ 
for an integer $l$ such that $0\leq l\leq m-1$ by Lemma \ref{T0305}\,(1). 
If $0\leq l\leq m-1$, we have that $3(3m-l-2)-2(3m-2)=3m-3l-2>0$, 
$2(3m-2)-3(l+m-1)=3m-3l-1>0$ and that 
\begin{align*}
\left| 1-\frac{3(3m-l-2)}{2(3m-2)}\right| & 
=\frac{3(3m-l-2)}{2(3m-2)}-1=\frac{3m-3l-2}{2(3m-2)} \\
{} & 
\leq \frac{3m-3l-1}{2(3m-2)}
=1-\frac{3(l+m-1)}{2(3m-2)}
=\left| 1-\frac{3(l+m-1)}{2(3m-2)}\right| . 
\end{align*}
Hence, 
\begin{align*}
\sum _{j=1}^n2^j\binom{n}{j}\left| 1-\frac{3j}{2n}\right| ^{2k} & 
\leq \sum _{j=1}^{m-2}a_{n, j}\left| 1-\frac{3j}{2n}\right| ^{2k} \\
{} & \quad 
+\sum _{l=0}^{m-1}(a_{3m-2, l+m-1}+a_{3m-2, 3m-l-2})
\left| 1-\frac{3(l+m-1)}{2(3m-2)}\right| ^{2k} \\
{} & 
\leq 10\sum _{j=1}^{2m-2}a_{n, j}\left| 1-\frac{3j}{2n}\right| ^{2k}
\leq 10\sum _{j=1}^{2m-2}\frac{2^jn^j}{j!}e^{-\frac{3jk}{n}} \\
{} & 
\leq 10\sum _{j=1}^{\infty}\frac{1}{j!}
e^{j\left( \log 2n-\frac{3k}{n}\right)}
\end{align*}
Finally, we assume $n=3m-3$ for an integer $m$ such that $m\geq 2$. 
For $l\in \{0, 1, \dots , m-1\}$, we see that $\disp{\frac{a_{3m-3, 3m-l-3}}
{a_{3m-3, l+m-1}}\leq 9}$ by Lemma \ref{T0305}\,(1). 
Moreover, for any integer $j$ with $0\leq l\leq m-2$, 
we have $(3m-l-3)-2(m-1)=m-l-1>0$, 
$2(m-1)-(l+m+1)=m-l-1>0$ and 
\begin{align*}
\left| 1-\frac{3(3m-l-3)}{2(3m-3)}\right| & 
=\frac{3m-l-3}{2(m-1)}-1=\frac{m-l-1}{2(m-1)} \\
{} & 
=1-\frac{l+m-1}{2(m-1)}
=\left| 1-\frac{3(l+m-1)}{2(3m-3)}\right| . 
\end{align*}
Hence, 
\begin{align*}
\sum _{j=1}^n2^j\binom{n}{j}\left| 1-\frac{3j}{2n}\right| ^{2k} & 
=\sum _{j=1}^{m-2}a_{n, j}\left| 1-\frac{3j}{2n}\right| ^{2k}
+a_{n, 2m-2}\left| 1-\frac{3(2m-2)}{2(3m-3)}\right| ^{2k} \\
{} & \quad 
+\sum _{l=0}^{m-2}(a_{3m-3, l+m-1}+a_{3m-3, 3m-l-3})
\left| 1-\frac{3(l+m-1)}{2(3m-3)}\right| ^{2k} \\
{} & 
\leq 10\sum _{j=1}^{2m-3}a_{n, j}\left| 1-\frac{3j}{2n}\right| ^{2k}
\leq 10\sum _{j=1}^{2m-3}\frac{2^jn^j}{j!}e^{-\frac{3jk}{n}} \\
{} & 
\leq 10\sum _{j=1}^{\infty}\frac{1}{j!}
e^{j\left( \log 2n-\frac{3k}{n}\right)}. 
\end{align*}
Therefore, for any integer $n$ such that $n\geq 3$, we have 
\begin{equation*}
\Vert \nu _n^{*k}-\pi _n\Vert _{TV}^2
\leq \frac{1}{4}\cdot 10\cdot \sum _{j=0}^{\infty}\frac{1}{j!}
e^{j\left( \log 2n-\frac{3k}{n}\right)}
=\frac{5}{2}\sum _{j=1}^{\infty}\frac{e^{-cj}}{j!}
=\frac{5}{2}(e^{e^{-c}}-1). 
\end{equation*}
\qed

\begin{theorem} \label{T0307}
Suppose $q=4$. 
For a positive integer $\disp{k=\frac{3n}{8}(\log 3n+c)}$ 
with $n\geq 2$ and $c>0$, we have 
\begin{equation} \label{0307}
\Vert \nu _n^{*k}-\pi _n\Vert _{TV}^2
\leq \frac{9}{4}(e^{e^{-c}}-1). 
\end{equation}
\end{theorem}

{\it{Proof}}.\,\,
Similarly to Theorem \ref{T0306}, we have 
\begin{align*}
\Vert \nu _n^{*k}-\pi _n\Vert _{TV}^2 & 
\leq \frac{1}{4}\sum _{j=1}^nd_j|{\cal{F}}(\nu _n)(\phi _j)|^{2k}
=\frac{1}{4}\sum _{j=1}^n3^j\binom{n}{j}
\left| 1-\frac{4j}{3n}\right| ^{2k}. 
\end{align*}
Put $\disp{a_{n, j}=3^j\binom{n}{j}}$. 
Assume $n=2m-1$ for some integer $m$ with $m\geq 2$. 
By Lemma \ref{T0305}\,(2), for an integer $j$ such that 
$\disp{0\leq l\leq \left\lfloor \frac{m-1}{2}\right\rfloor}$, we have 
$\disp{\frac{a_{2m-1, 2m-l-1}}{a_{2m-1, l+m-1}}\leq 8}$. 
Moreover, if $\disp{0\leq l\leq \left\lfloor \frac{m-1}{2}\right\rfloor}$, we see that 
\begin{align*}
4(2m-l-1) & 
\geq 4\left( 2m-\left( \frac{m-1}{2}\right) -1\right) =6m-2\geq 3(2m-1) \\
{} & 
\geq 6m-6=4\left( \frac{m-1}{2}+m-1\right) \geq 4(l+m-1). 
\end{align*}
So we have 
\begin{align*}
\left| 1-\frac{4(2m-l-1)}{3(2m-1)}\right| & 
=\frac{4(2m-l-1)}{3(2m-1)}-1=\frac{2m-4l-1}{3(2m-1)} \\
{} & 
\leq \frac{2m-4l+1}{3(2m-1)}=1-\frac{4(l+m-1)}{3(2m-1)}
=\left| 1-\frac{4(l+m-1)}{3(2m-1)}\right| . 
\end{align*}
Hence, if $m=2r$ for some positive integer $r$, we have 
\begin{align*}
\sum _{j=1}^n3^j\binom{n}{j}\left| 1-\frac{4j}{3n}\right| ^{2k} & 
\leq \sum _{j=1}^{m-2}a_{n, j}\left| 1-\frac{4j}{3n}\right| ^{2k}
+a_{n, 3r-1}\left| 1-\frac{4(3r-1)}{3n}\right| ^{2k} \\
{} & \quad 
+\sum _{l=0}^{r-1}(a_{2m-1, l+m-1}+a_{2m-1, 2m-l-1})
\left| 1-\frac{4(l+m-1)}{3(2m-1)}\right| ^{2k} \\
{} & 
\leq 9\sum _{j=1}^{3r-1}a_{n, j}\left| 1-\frac{4j}{3n}\right| ^{2k}
\leq 9\sum _{j=1}^{3r-1}\frac{3^jn^j}{j!}e^{-\frac{8jk}{3n}}
\leq 9\sum _{j=1}^{\infty}\frac{1}{j!}
e^{j\left( \log 3n-\frac{8k}{3n}\right)}.   
\end{align*}
If $m=2r+1$ for a positive integer $r$, we have $\disp{r=\frac{m-1}{2}}$ and 
\begin{align*}
\sum _{j=1}^n3^j\binom{n}{j}\left| 1-\frac{4j}{3n}\right| ^{2k} & 
\leq \sum _{j=1}^{m-2}a_{n, j}\left| 1-\frac{4j}{3n}\right| ^{2k} \\
{} & \quad 
+\sum _{l=0}^r(a_{2m-1, l+m-1}+a_{2m-1, 2m-l-1})
\left| 1-\frac{4(l+m-1)}{3(2m-1)}\right| ^{2k} \\
{} & 
\leq 9\sum _{j=1}^{3r}a_{n, j}\left| 1-\frac{4j}{3n}\right| ^{2k}
\leq 9\sum _{j=1}^{3r}\frac{3^jn^j}{j!}e^{-\frac{8jk}{3n}}
\leq 9\sum _{j=1}^{\infty}\frac{1}{j!}
e^{j\left( \log 3n-\frac{8k}{3n}\right)}.  
\end{align*}
Next, we consider the case $n=2m-2$ for some integer $m$ with $m\geq 2$. 
For an integer $l$ such that $\disp{0\leq l\leq \left\lfloor 
\frac{m-2}{2}\right\rfloor}$, we have 
$\disp{\frac{a_{2m-2, 2m-l-2}}{a_{2m-2, l+m-1}}\leq 8}$ by Lemma 
\ref{T0305}\,(2), and that 
\begin{align*}
2(2m-l-2) & 
\geq 2\left( 2m-\left( \frac{m-2}{2}\right) -2\right)=3m-2\geq 3(m-1) \\
{} & 
\geq 3m-4=2\left( \left( \frac{m-2}{2}\right) +m-1\right) \geq 2(l+m-1).  
\end{align*}
Hence, 
\begin{align*}
\left| 1-\frac{4(2m-l-2)}{3(2m-2)}\right| & 
=\frac{2(2m-l-2)}{3m-3}-1=\frac{m-2l-1}{3m-3} \\
{} & 
=1-\frac{2(l+m-1)}{3m-3}=\left| 1-\frac{4(l+m-1)}{3(2m-2)}\right| . 
\end{align*}
If $m=2r$ for some positive integer $r$, we have $n=4r-2$ and 
\begin{align*}
\sum _{j=1}^n3^j\binom{n}{j}\left| 1-\frac{4j}{3n}\right| ^{2k} & 
=\sum _{j=1}^{m-2}a_{n, j}\left| 1-\frac{4j}{3n}\right| ^{2k} \\
{} & \quad 
+\sum _{l=0}^{r-1}(a_{2m-2, l+m-1}+a_{2m-2, 2m-l-2})
\left| 1-\frac{4(l+m-1)}{3(2m-2)}\right| ^{2k} \\
{} & 
\leq 9\sum _{j=1}^{3r-2}a_{n, j}\left| 1-\frac{4j}{3n}\right| ^{2k}
\leq 9\sum _{j=1}^{3r-2}\frac{3^jn^j}{j!}e^{-\frac{8jk}{3n}}
\leq 9\sum _{j=1}^{\infty}\frac{1}{j!}
e^{j\left( \log 3n-\frac{8k}{3n}\right)}.  
\end{align*}
If there exists a positive integer $r$ such that $m=2r+1$, 
then $n=4r$, $\disp{r=\frac{m-1}{2}}$ and 
\begin{align*}
\sum _{j=1}^n3^j\binom{n}{j}\left| 1-\frac{4j}{3n}\right| ^{2k} & 
=\sum _{j=1}^{m-2}a_{n, j}\left| 1-\frac{4j}{3n}\right| ^{2k}
+a_{2m-2, 3r}\left| 1-\frac{4\cdot 3r}{3\cdot 4r}\right| ^{2k} \\
{} & \quad 
+\sum _{l=0}^{r-1}(a_{2m-2, l+m-1}+a_{2m-2, 2m-l-2})
\left| 1-\frac{4(l+m-1)}{3(2m-2)}\right| ^{2k} \\
{} & 
\leq 9\sum _{j=1}^{3r-1}a_{n, j}\left| 1-\frac{4j}{3n}\right| ^{2k}
\leq 9\sum _{j=1}^{3r-1}\frac{3^jn^j}{j!}e^{-\frac{8jk}{3n}}
\leq 9\sum _{j=1}^{\infty}\frac{1}{j!}
e^{j\left( \log 3n-\frac{8k}{3n}\right)}.  
\end{align*}
Therefore, for any integer $n$ with $n\geq 2$, we have 
\begin{equation*}
\Vert \nu _n^{*k}-\pi _n\Vert _{TV}^2
\leq \frac{1}{4}\cdot 9\sum _{j=1}^{\infty}\frac{1}{j!}
e^{j\left( \log 3n-\frac{8k}{3n}\right)}
=\frac{9}{4}\sum _{j=1}^{\infty}\frac{e^{-cj}}{j!}
=\frac{9}{4}(e^{e^{-c}}-1). 
\end{equation*}
\qed

\begin{remark}
{\rm{Diaconis and Ram study the analysis on the deformation of Markov chains 
on Coxeter groups, called the {\it{Metropolis algorithm}}. 
The random walks on the Hamming schemes $\{H(n, q)\}_{n=1}^{\infty}$ are 
the special cases of them on the hypercubes $\{H(n, 2)\}_{n=1}^{\infty}$ 
with $\disp{\theta =\frac{1}{q-1}}$ (\cite{DR}, Theorem 5.4). 
However, our majorant functions for the total variance distances are 
simpler than that they give.}} 
\end{remark}

\section{Lower bounds.}
\quad\,\,
In this section we give a minorant function for the total variance distance 
$\Vert \nu _n^{*k}-\pi _n\Vert _{TV}$. 
Let $\mu$ be a probabilitiy measure on $X_n$ and 
$f$ a function on $X_n$. 
We denote by $E_{\mu}(f)$ the {\it{expectation}} of $f$ 
with respect to $\mu$, that is, 
\begin{equation*}
E_{\mu}(f)=\sum _{x\in X_n}f(x)\mu (x). 
\end{equation*}
The {\it{variance}} $Var_{\mu}(f)$ of $f$ respect to $\mu$ is defined by 
\begin{equation*}
Var _{\mu}(f)=E_{\mu}((f-E_{\mu}(f))^2)
=E_{\mu}(f^2)-E_{\mu}(f)^2. 
\end{equation*}
In order to compute expectations and variances of spherical functions, 
we need the precise description of $\phi _j$'s. 
We see that $\phi _0=1$, the constant function with the value $1$. 
The spherical function $\phi _1$ is described as 
\begin{align*}
\phi _1(l) & 
=1+\frac{(-1)\cdot (-l)}{(-n)\cdot 1}\cdot \frac{q}{q-1}
=1-\frac{lq}{n(q-1)}, 
\end{align*}
where $l\in \{0, 1, \dots , n\}$. 
We calculate the spherical function $\phi _2$ to compute 
the variance of $\phi _1$. 
For $l\in \{0, 1, \dots , n\}$, we have 
\begin{align*}
\phi _2(l) & 
=1+\frac{(-2)\cdot (-l)}{(-n)\cdot 1}\cdot \frac{q}{q-1}
+\frac{(-2)_2(-l)_2}{(-n)_2\cdot 2!}\cdot \left( \frac{q}{q-1}\right) ^2 \\
{} & 
=1-\frac{2lq}{n(q-1)}+\frac{l(l-1)q^2}{n(n-1)(q-1)^2} \\
{} & 
=1-\frac{lq((2n-1)q-(2n-2))}{n(n-1)(q-1)^2}+\frac{l^2q^2}{n(n-1)(q-1)^2}.  
\end{align*}

\begin{lemma} \label{T0401}
We have 
\begin{equation} \label{0401}
\phi _1^2=\frac{1}{n(q-1)}\phi _0+\frac{q-2}{n(q-1)}\phi _1
+\frac{n-1}{n}\phi _2. 
\end{equation}
\end{lemma}

{\it{Proof}}.\,\,
For $l\in \{1, \dots , n\}$, we see that  
\begin{equation*} 
\phi _1(l)^2=\left( 1-\frac{lq}{n(q-1)}\right) ^2
=1-\frac{2lq}{n(q-1)}+\frac{l^2q^2}{n^2(q-1)^2}. 
\end{equation*}
On the other hand, we have 
\begin{align*}
{} & \frac{1}{n(q-1)}\phi _0(l)+\frac{q-2}{n(q-1)}\phi _1(l)
+\frac{n-1}{n}\phi _2(l) \\
{} & \quad 
=\frac{1}{n(q-1)}+\frac{q-2}{n(q-1)}\left( 1-\frac{lq}{n(q-1)}\right) \\
{} & \quad\quad
+\frac{n-1}{n}\left( 1-\frac{lq((2n-1)q-(2n-2))}{n(n-1)(q-1)^2}
+\frac{l^2q^2}{n(n-1)(q-1)^2}\right) \\
{} & \quad 
=\frac{1+(q-2)+(n-1)(q-1)}{n(q-1)}
-\frac{lq((q-2)+((2n-1)q-(2n-2)))}{n^2(q-1)^2}
+\frac{l^2q^2}{n^2(q-1)^2} \\
{} & \quad 
=1-\frac{2lq}{n(q-1)}+\frac{l^2q^2}{n^2(q-1)^2}, 
\end{align*}
it completes the proof of the lemma. 
\qed

\vspace{1pc}

Here, we compute the expectations of $\phi _j$'s and 
the variance of $\phi _1$ with respect to $\pi _n$. 

\begin{lemma} \label{T0402}
{\rm{(1)}}\,\,
We have 
\begin{equation} \label{0402}
E_{\pi _n}(\phi _j)=\begin{cases}
1, & j=1,           \\
0, & 1\leq j\leq n. 
\end{cases}
\end{equation}
{\rm{(2)}}\,\,
One has $\disp{Var _{\pi _n}(\phi _1)=E_{\pi _n}(\phi _1^2)=\frac{1}{n(q-1)}}$. 
\end{lemma}

{\it{Proof}}.\,\,
(1)\,\,
We see that 
\begin{equation*}
E_{\pi _n}(\phi _0)=\sum _{x\in X_n}\phi _0(x)\pi _n(x)
=(\sharp X_n)\cdot 1\cdot \frac{1}{\sharp X_n}=1. 
\end{equation*}
We compute $E_{\pi _n}(\phi _j)$ for $j\in \{1, 2, \dots , n\}$. 
For each $x\in X_n$, there exists $l\in \{0, 1, \dots , n\}$ such that 
$x\in H_n\cdot x^{(l)}$. 
We see that 
\begin{equation*}
\sharp (H_n\cdot x^{(l)})
=\sharp \{x\in X_n\,;\,d(x^{(0)}, x)=l\}=(q-1)^l\binom{n}{l}. 
\end{equation*}
Since $\sharp X_n=q^n$, we have 
\begin{align*}
(\sharp X_n)E_{\pi _n}(\phi _j) & 
=(\sharp X_n)\sum _{x\in X_n}\phi _j(x)\pi _n(x)
=(\sharp X_n)\cdot \frac{1}{\sharp X_n}\sum _{l=0}^n
\sum _{x\in H_n\cdot x^{(l)}}\phi _j(x) \\
{} & 
=\sum _{l=0}^n(q-1)^l\binom{n}{l}\sum _{r=0}^j
\frac{(-j)_r(-l)_r}{(-n)_rr!}\left( \frac{q}{q-1}\right) ^r \\
{} & 
=\sum _{r=0}^j\frac{(-j)_r}{(-n)_rr!}\left( \frac{q}{q-1}\right) ^r
\sum _{l=0}^n(q-1)^l(-l)_r\binom{n}{l}. 
\end{align*}
For $r\in \{0, 1, \dots , j\}$, we see that 
\begin{align*}
\sum _{l=0}^n(q-1)^l(-l)_r\binom{n}{l} & 
=(-1)^rn(n-1)\cdots (n-r+1)(q-1)^r\sum _{l=r}^n(q-1)^{l-r}\binom{n-r}{l-r} \\
{} & 
=(-n)_r(q-1)^r\sum _{l=0}^{n-r}(q-1)^l\binom{n-r}{l}
=(-n)_rq^{n-r}(q-1)^r. 
\end{align*}
Hence, we have 
\begin{align*}
(\sharp X_n)E_{\pi _n}(\phi _j) & 
=\sum _{r=0}^j\frac{(-j)_r}{(-n)_rr!}\left( \frac{q}{q-1}\right) ^r
\cdot (-n)_rq^{n-r}(q-1)^r \\
{} & 
=q^n\sum _{r=0}^j\frac{(-j)_r}{r!}
=q^n\sum _{r=0}^j(-1)^r\binom{j}{r}=0. 
\end{align*}
(2)\,\,
By Lemma \ref{T0401}, we have 
\begin{align*}
Var _{\pi _n}(\phi _1) & 
=E_{\pi _n}(\phi _1^2)-E_{\pi _n}(\phi _1)^2=E_{\pi _n}(\phi _1^2) \\
{} & 
=\frac{1}{n(q-1)}E_{\pi _n}(\phi _0)-\frac{q-2}{n(q-1)}E_{\pi _n}(\phi _1)
+\frac{n-1}{n}E_{\pi _n}(\phi _2) \\
{} & 
=\frac{1}{n(q-1)}.
\end{align*}
\qed

\vspace{1pc}

Next, we calculate the expectations of $\phi _j$'s and 
estimate the variance of $\phi _1$ with respect to $\nu _n^{*k}$ 
for any nonnegative integer $k$. 

\begin{lemma} \label{T0403}
Let $n$ be a positive integer, $q$ an integer with $q\geq 2$ and 
$k$ a nonnegative integer. 
\newline
{\rm{(1)}}
For $j\in \{0, 1, \dots , n\}$, we have 
\begin{equation} \label{0403}
E_{\nu _n^{*k}}(\phi _j)=\left( 1-\frac{jq}{n(q-1)}\right) ^k. 
\end{equation}
{\rm{(2)}}\,\,
Assume that $(n-2)(q-1)\geq 2$. 
Then, one has 
\begin{equation} \label{0404}
Var _{\nu _n^{*k}}(\phi _1)\leq \frac{1}{n}. 
\end{equation}
\end{lemma}

{\it{Proof}}.\,\,
(1)\,\,
Since $\phi _j$ is real-valued for any $j\in \{0, 1, \dots , n\}$, we have 
\begin{align*}
E_{\nu _n^{*k}}(\phi _j) & 
=\sum _{x\in X_n}\phi _j(x)\nu _n^{*k}(x)
=\sum _{x\in X_n}\nu _n^{*k}(x)\overline{\phi _j(x)} \\
{} & 
={\cal{F}}(\nu _n^{*k})(\phi _j)
={\cal{F}}(\nu _n)(\phi _j)^k
=\left( 1-\frac{jq}{n(q-1)}\right) ^k. 
\end{align*}
(2)\,\,
By Lemma \ref{T0401}, we see that 
\begin{align*}
Var _{\nu _n^{*k}}(\phi _1) & 
=E_{\nu _n^{*k}}(\phi _1^2)-E_{\nu _n^{*k}}(\phi _1)^2 \\
{} & 
=\frac{1}{n(q-1)}E_{\nu _n^{*k}}(\phi _0)
+\frac{q-2}{n(q-1)}E_{\nu _n^{*k}}(\phi _1)
+\frac{n-1}{n}E_{\nu _n^{*k}}(\phi _2)-E_{\nu _n^{*k}}(\phi _1)^2 \\
{} & 
=\frac{1}{n(q-1)}+\frac{q-2}{n(q-1)}\left( 1-\frac{q}{n(q-1)}\right) ^k
+\frac{n-1}{n}\left( 1-\frac{2q}{n(q-1)}\right) ^k \\
{} & \quad 
-\left( 1-\frac{q}{n(q-1)}\right) ^{2k}. 
\end{align*}

Since $n(q-1)-2q=(n-2)(q-1)-2\geq 0$, we see that 
\begin{equation*}
0\leq 1-\frac{2q}{n(q-1)}\leq \left( 1-\frac{q}{n(q-1)}\right) ^2.
\end{equation*} 
Hence, we have 
\begin{align*}
Var _{\nu _n^{*k}}(\phi _1) & 
\leq \frac{1}{n(q-1)}+\frac{q-2}{n(q-1)}\left( 1-\frac{q}{n(q-1)}\right) ^k
-\frac{1}{n}\left( 1-\frac{q}{n(q-1)}\right) ^{2k} \\
{} & 
\leq \frac{1}{n(q-1)}+\frac{q-2}{n(q-1)}\left( 1-\frac{q}{n(q-1)}\right) ^k
\leq \frac{1}{n(q-1)}+\frac{q-2}{n(q-1)}=\frac{1}{n}, 
\end{align*}
since $n(q-1)-q>n(q-1)-2q\geq 0$. 
\qed

\vspace{1pc}

Now, we give a minorant function for total variance distance. 

\begin{theorem} \label{T0404}
Assume that $q\geq 2$.
We fix a positive real number $c_0>0$. 
For any positive real number $b>0$, there exists a positive integer 
$n_0$ such that for any integer $n$ with $n\geq n_0$ and any integer 
$\disp{k=\frac{n(q-1)}{2q}(\log n(q-1)-c)}$ 
with $0\leq c\leq \min \{c_0, \log n(q-1)\}$, we have 
\begin{equation} \label{0405}
\Vert \nu _n^{*k}-\pi _n\Vert _{TV}\geq 1-(4q+b)e^{-c}. 
\end{equation}
\end{theorem}

{\it{Proof}}.\,\,
We write $\disp{n_1=\max \left\{ \left\lceil \frac{e^{c_0}}{q-1}\right\rceil , 
4\right\}}$ for simplicity, and assume that $n\geq n_1\geq 4$. 
Then, we see that $c_0\leq \log n_1(q-1)\leq \log n(q-1)$. 
Put 
\begin{equation} \label{0406}
\beta _{n, k}=\sqrt{\frac{q}{(4q+b)(q-1)}}e^{\frac{c}{2}}. 
\end{equation}
Using $\beta _{n, k}$, we define a subset 
$B_{n, k}\subset X_n$ of $X_n$ by 
\begin{equation} \label{0407}
B_{n, k}=\left\{ x\in X_n\,;\,|\phi _1(x)|
<\frac{\beta _{n, k}}{\sqrt{n}}\right\} . 
\end{equation}
By Markov's inequality, we have 
\begin{align*}
\pi _n(B_{n, k}) & 
=1-\pi _n\left( \left\{x\in X_n\,;\,|\phi _1(x)|
\geq \frac{\beta _{n, k}}{\sqrt{n}}\right\} \right) \\
{} & 
=1-\pi _n\left( \left\{ x\in X_n\,;\,\phi _1(x)^2
\geq \frac{\beta _{n, k}^2}{n}\right\} \right) \\
{} & 
\geq 1-\frac{n}{\beta _{n, k}^2}E_{\pi _n}(\phi _1^2)
=1-\frac{n}{\beta _{n, k}^2}\cdot \frac{1}{n(q-1)}
=1-\frac{1}{\beta _{n, k}^2(q-1)}. 
\end{align*}
We define a function $\omega\,:\,[0, 1)\longrightarrow {\mathbb{R}}$ 
on the interval $[0, 1)$ by 
\begin{equation*}
\log (1-x)=-x-\frac{x^2}{2}\omega (x), 
\end{equation*}
where $x\in [0, 1)$.
Then, $\omega (x)\geq 0$ for any $x\in [0, 1)$, $\disp{\lim _{x\to 0}
\omega (x)=1}$ and we have \begin{align*}
E_{\nu _n^{*k}}(\phi _1) & 
=\left( 1-\frac{q}{n(q-1)}\right) ^k
=\exp \left( \log \left( 1-\frac{q}{n(q-1)}\right) ^k\right) \\
{} & 
=\exp \left( \log \left( 1-\frac{q}{n(q-1)}\right) 
\cdot \frac{n(q-1)}{2q}(\log n(q-1)-c)\right) \\
{} &  
=\exp \left( \left( -\frac{q}{n(q-1)}-\frac{q^2}{2n^2(q-1)^2}
\omega \left( \frac{q}{n(q-1)}\right) \right) 
\frac{n(q-1)}{2q}(\log n(q-1)-c)\right) \\
{} & 
=\frac{e^{\frac{c}{2}}}{\sqrt{n(q-1)}}
\exp \left( \frac{q(c-\log n(q-1))}{4n(q-1)}
\omega \left( \frac{q}{n(q-1)}\right) \right) . 
\end{align*}
Since $0\leq c\leq c_0\leq \log n_1(q-1)\leq \log n(q-1)$ and 
$\disp{\lim _{n\to \infty}
\frac{\log n(q-1)}{n}=\lim _{n\to \infty}\frac{1}{n}=0}$, 
we see that 
\begin{align*}
{} & 
0<\exp \left( -\frac{q\log n(q-1)}{4n(q-1)}
\omega \left( \frac{q}{n(q-1)}\right) \right) \\
{} & \,\,\,
\leq \exp \left( \frac{q(c-\log n(q-1))}{4n(q-1)}
\omega \left( \frac{q}{n(q-1)}\right) \right) \\
{} & \,\,\,
\leq \exp \left( \frac{q(c_0-\log n(q-1))}{4(q-1)}
\omega \left( \frac{q}{n(q-1)}\right) \right) \leq 1, \\
{} & 
\lim _{n\to \infty}\exp \left( \frac{-q\log n(q-1)}{4n(q-1)}
\omega \left( \frac{q}{n(q-1)}\right) \right) =1. 
\end{align*}
Hence, there exists an integer $n_0$ with $n_0\geq n_1$ such that 
for any integer $n$ with $n\geq n_0$, we have 
\begin{equation*}
\exp \left( \frac{-q\log n(q-1)}{4n(q-1)}
\omega \left( \frac{q}{n(q-1)}\right) \right) 
\geq 2\sqrt{\frac{q}{4q+b}},  
\end{equation*}
it implies that 
\begin{equation*}
E_{\nu _n^{*k}}(\phi _1)\geq 2\sqrt{\frac{q}{n(4q+b)(q-1)}}
e^{\frac{c}{2}}=\frac{2\beta _{n, k}}{\sqrt{n}}. 
\end{equation*}
So we have 
\begin{equation*}
B_{n, k}\subset B'_{n, k}=\left\{ x\in X_n\,;\,
|\phi _1(x)-E_{\nu _n^{*k}}(\phi _1)|
\geq E_{\nu _n^{*k}}(\phi _1)-\frac{\beta _{n, k}}{\sqrt{n}}\right\} . 
\end{equation*}
Hence, by Chebyshev's inequality, 
\begin{equation*}
\nu _n^{*k}(B_{n, k})\leq \nu _n^{*k}(B'_{n, k})
\leq \frac{Var _{\nu _n^{*k}}(\phi _1)}
{\disp{\left( E_{\nu _n^{*k}}(\phi _1)-\frac{\beta _{n, k}}{\sqrt{n}}\right) ^2}}
\leq \frac{\disp{\,\,\,\frac{1}{n}\,\,\,}}
{\disp{\,\,\frac{\beta _{n, k}^2}{n}\,\,}}
=\frac{1}{\beta _{n, k}^2}. 
\end{equation*}
Therefore
\begin{align*}
\Vert \nu _n^{*k}-\pi _n\Vert _{TV} & 
\geq \pi _n(B_{n, k})-\nu _n^{*k}(B_{n, k})
\geq 1-\frac{1}{\beta _{n, k}^2(q-1)}-\frac{1}{\beta _{n, k}^2} \\
{} & 
=1-\frac{q}{q-1}\cdot \frac{(4q+b)(q-1)}{q}e^{-c}=1-(4q+b)e^{-c}. 
\end{align*}
\qed

\vspace{1pc}

The above theorem gives a simple lower bound 
of $\Vert \nu _n^{*k}-\pi _n\Vert _{TV}$. 

\begin{corollary} \label{T0405}
Fix a positive real number $c>0$. 
Put $\disp{a_n=\frac{n(q-1)}{2q}\log n(q-1)}$ and 
$\disp{b_n=\frac{n(q-1)}{2q}}$. 
Then, we have 
\begin{equation} \label{0408}
\liminf _{n\to \infty}\Vert \nu _n^{*\lfloor a_n-cb_n\rfloor}
-\pi _n\Vert _{TV}\geq 1-4qe^{-c}. 
\end{equation}
\end{corollary}

\begin{remark}
{\rm{In Theorem 5.8 of \cite{M2}, Mizukawa gives a minorant function 
for random walks with staying on $(K/L)^n$,  
where $(K, L)$ is a Gelfand pair. 
Our object is the case where $(K, L)=(S_q, S_{q-1})$, 
$\disp{a_0=\frac{1}{q-1}}$ and $mp=1$. 
Put $\disp{\gamma =2\sqrt{\frac{q}{4q+b}}}$. 
Then we have that $0<\gamma <1$ and that $\disp{\delta 
=\frac{4(q-1)(a_0+1)}{\gamma ^2}}$ $=4q+b$}}. 

\end{remark}

\vspace{1pc}

\begin{flushright}
Katsuhiko Kikuchi \\
Department of Mathematics \\
Kyoto University \\
606-8502 Kyoto, JAPAN \\
e-mail\,:\,kikuchi@math.kyoto-u.ac.jp
\end{flushright}

\end{document}